\input amstex
    \documentstyle{amsppt}
    \magnification=1200
    \pagewidth{13cm}
    \voffset=-1in
    \topmatter
    \title
    Bounds for the cup-length of Poincar\'e spaces
    and their applications\endtitle
    \author
    J\'ulius Korba\v s
    \endauthor
    \address
    Department of Algebra, Geometry, and Mathematical Education,
    Faculty of Mathematics, Physics, and Informatics,
    Comenius University,
    Mlynsk{\'a} dolina,
    SK-842 48 Bratislava 4,
    Slovakia
    \endaddress
    \email korbas\@fmph.uniba.sk \newline
    {\it or}\endemail
    \address
    Mathematical Institute, Slovak Academy of Sciences,
    {\v S}tef{\'a}nikova 49, SK-814 73 Bratislava 1, Slovakia \endaddress
    \email matekorb\@savba.sk \endemail
    \abstract Our main result offers a new (quite systematic) way of
    deriving bounds for the cup-length of Poincar\'e spaces
    over fields; we outline a general research program based on this result.
    For the oriented Grassmann manifolds, already a limited
    realization of the program leads, in many cases, to the
    exact values of the cup-length and to interesting information on the
    Lyusternik-Shnirel'man category.

    \endabstract
    \thanks The author was supported in part by two grants
    of VEGA (Slovakia)
    \endthanks
    \leftheadtext{J\'ulius Korba\v s}
    \rightheadtext{Bounds for the cup-length and their applications}
     \endtopmatter
    \document
    \TagsOnRight

\head 1. INTRODUCTION AND STATEMENT OF RESULTS\endhead

The cup-length over a field $R$, $\operatorname{cup}_{R}(X)$, of a path connected topological
space $X$ is the supremum of all numbers $c$ such that there exist positive dimensional
cohomology classes $a_1,\dots,a_c\in H^\ast(X; R)$ such that their cup product
$a_1\cup\dots\cup a_c$ is nonzero.
Instead of the standard notation $a\cup b$, we shall
mostly write $a\cdot b$ or just $ab$.
Although the main topic of this paper is the cup-length, we
shall also keep in mind its relation to another important
invariant: we have $\operatorname{cat}(X)\geq
1+\operatorname{cup}_R(X)$, where $\operatorname{cat}(X)$ is the L-S category (the
Lyusternik-Shnirel'man category). We recall that $\operatorname{cat}(X)$ is defined to be the
least positive integer $k$ such that $X$ can be covered by $k$ open subsets each of which is
contractible in $X$. If no such integer exists, then one puts $\operatorname{cat}(X)=\infty$.
Note that some authors (see e.g. \cite{5}) prefer to modify the definition by subtracting $1$
from our value of $\operatorname{cat}(X)$.

The L-S category is very hard to compute; a longstanding problem in topology (cf. Ganea's
list \cite{7}) is the task to find its value for familiar manifolds. In general, the
cup-length is difficult to calculate, and (a formula for) its value remains unknown for many
commonly used spaces, e.g. for the great majority of the Grassmann manifolds $O(n)/O(k)\times
O(n-k)$; explicit formulae for $\operatorname{cup}_{\Bbb Z_2}(O(n)/O(k)\times O(n-k))$ are
known for $k\leq 5$, with some exceptions (cf. R. Stong's \cite{12}).

\subhead 1.1 Statement of the main result \endsubhead When trying to find the cup-length of a
space, one might start from calculations of the height for some nonzero elements in its
cohomology. Recall that if $x\in H^i(X; R)$ ($i>0$) is a nonzero cohomology class, then the
height of $x$ over $R$, denoted by $\operatorname{ht}_R(x)$, is the supremum of all the
numbers $k$ such that $x^k\neq 0$; of course, one has $\operatorname{cup}_R(X)\geq
\operatorname{ht}_R(x)$. Our main result is the following general theorem which we shall
prove in Section 2. Its most interesting part, (b), shows that if $X$ is an $R$-Poincar\'e
space, then height-related information may lead also to {\it upper} bounds for
$\operatorname{cup}_R(X)$; part (a) (in particular (a2) and (a3)) is almost obvious, but we
include it in order to have the theorem as a comfortable reference tool.

\proclaim{Theorem A} Let $R$ be a field, and $X\neq\emptyset$ be a path-connected
$R$-Poincar\'e space of formal dimension $n$. Let the first two nonzero reduced
$R$-cohomology groups of $X$, in dimensions less than $n$, be $\tilde H^r(X; R)$ and $\tilde
H^q(X; R)$, $r\leq q<n$. Then the cup-length satisfies the following. \roster
\item"(a)" One always has
$$\operatorname{cup}_R(X)\leq \frac{n}{r}. \tag a1 $$
If there is a cohomology class $x\in H^r(X; R)\setminus
\{0\}$ such that $r\cdot\operatorname{ht}_R(x)=n$, then
$$\operatorname{cup}_R(X)= \frac{n}{r}. \tag a2 $$
If one finds $a_1,\dots,a_s\in H^r(X; R)\setminus \{0\}$ such that $a_1^{t_1}\cdots
a_s^{t_s}\neq 0$ for some non-negative integers $t_1,\dots,t_s$ such that $t_1+\dots+t_s>0$
and $r(t_1+\dots +t_s)<n$, then one has
$$\operatorname{cup}_R(X)\geq 1+t_1+\dots +t_s. \tag a3 $$

\item"(b)"
Suppose that $r<q$ and there exists a basis
$a(1)_r,\dots, a(t)_r$ for the $R$-vector space
$H^r(X; R)$ such that
$$a(1)_r^{k_1+1}=0,\dots, a(t)_r^{k_t+1}=0$$
for some positive integers $k_1,\dots, k_t$ satisfying the condition
$r(k_1+\dots+k_t)<n$.

Then the upper bound given in $(\text{\rm a})$ improves to
$$\operatorname{cup}_R(X)\leq k_1+\dots+k_t+
\left[\frac{n-r(k_1+\dots+k_t)}{q}\right]< \frac{n}{r}. \tag b1 $$
\endroster
\endproclaim
We remark that if $X$ has the homotopy type of an $(r-1)$-connected ($r\geq 1$) finite
CW-complex, then, as proved by Grossman in \cite{8}, $\operatorname{cat}(X)\leq
1+\dsize\frac{\operatorname{dim}(X)}{r}$; as a consequence, one then also has
$\operatorname{cup}_R(X)\leq \dsize\frac{\operatorname{dim}(X)}{r}$. Note that Grossman's
inequality and our (a1) in Theorem A coincide for $X$ having the homotopy type of an
$(r-1)$-connected ($r\geq 1$) finite CW-complex if $X$ is an $R$-Poincar\'e space; but our
(a1) requires just $R$-homological $(r-1)$-connectedness, while for Grossman's upper bound
the standard $(r-1)$-connectedness is required.

Theorem A can serve as a basis for the following research program. \proclaim{Research
Program} Let $X$ be an $R$-Poincar\'e space such as we suppose in the theorem, such that
$r<q$ {\rm(}if $r=q$, then one readily sees that $\operatorname{cup}_R(X)=2${\rm)}, and such
that $H^r(X; R)$ is finitely generated as a vector space over $R$. Then \roster
\item"(1)" use $(\text{\rm a}1)$ to calculate the initial upper bound for $\operatorname{cup}_R(X)$;
\item"(2)" study vanishing of products of elements in $H^*(X; R)$ and find {\rm (}possibly using
$(\text{\rm a}2)$ or $(\text{\rm a}3)${\rm)} a lower bound, as high as you can, for
$\operatorname{cup}_R(X)$;
\item"(3)" if the best upper and lower bounds you have obtained do not coincide,
then find bases {\rm(}as many as you can{\rm)}
in $H^r(X; R)$, study vanishing of powers of their
elements and, when possible, apply $(\text{\rm b}1)$ to obtain a better upper bound.
\endroster
\endproclaim
Note that this Research Program can be applied, in particular, to any $R$-orientable closed
$n$-dimensional manifold with $R$-cohomology different from the $R$-cohomology of the
$n$-sphere $S^n$.

\subhead 1.2 Statement of some applications of Theorem A\endsubhead To illustrate its
usefulness, we shall use Theorem A, with $r=2$, $R=\Bbb Z_2$ or $r=4$, $R=\Bbb Q$, for
deriving several new results on the cup-length of the oriented Grassmann manifolds; at the
same time, we shall pay attention to their L-S category. In many cases, in spite of knowing
the cup-length, one might be very far from knowing the exact value of the L-S category (see
e.g. \cite{5}). But we shall show that for the oriented Grassmann manifolds $\tilde
G_{n,k}=SO(n)/SO(k)\times SO(n-k)$, with $n\geq 2k\geq 6$, the cup-length yields, at least in
some cases, a good amount of information about the L-S category (we shall take $k\geq 3$,
because for $\tilde G_{n,1}=S^{n-1}$ and for $\tilde G_{n,2}$ (complex quadrics) the
cup-length and L-S category are known). More precisely, already our limited (in other words:
just illustrative) realization of the Research Program brings new non-trivial estimates for
the cup-length; in particular, we find the exact values of $\operatorname{cup}_{\Bbb
Z_2}(\tilde G_{6,3})$ and also of $\operatorname{cup}_{\Bbb Q}(\tilde G_{n,k})$ for
infinitely many $(n,k)$ with $k\geq 4$. At the same time, it turns out that for $\tilde
G_{n,3}$ the cup-length lower bound (that is, $1$ plus the cup-length) and the L-S category
can be very close to each other. In two cases, for $\tilde G_{6,3}$ and $\tilde G_{7,3}$, our
cup-length lower bound and the Grossman upper bound differ by just $1$, for $\tilde G_{8,3}$
the difference is $2$.

Now we pass to detailed statements of the applications which we have roughly outlined. In
order to simplify the notation, we shall write $\operatorname{cup}(X)$ instead of
$\operatorname{cup}_{\Bbb Z_2}(X)$, and $\operatorname{ht}(x)$ instead of
$\operatorname{ht}_{\Bbb Z_2}(x)$.

Using Theorem A(a) (mainly (a3); cf. (2) in the Research Program), we prove the following new
lower bounds for the cup-length of the oriented Grassmann manifolds in Section~3.

\proclaim{Proposition B} For the oriented Grassmann manifolds $\tilde G_{n,k}$, $n\geq 2k\geq
6$, we have
 \roster
\item"(a)"$\operatorname{cup}(\tilde G_{6,3})\geq 3$;
\item"(b)" $\operatorname{cup}(\tilde G_{n,3})\geq \cases
\dsize\frac{n+3}{2} &\text{ if } n\geq 7 \text{ is odd, } n\notin \{9,11\}, \\
\dsize\frac{n+2}{2} &\text{ if } n\geq 8 \text{ is even, } n\notin
\{10,12\}. \\
\endcases $
\item"(c)" Each of $\operatorname{cup}(\tilde G_{9,3})$, $\operatorname{cup}(\tilde
G_{10,3})$, $\operatorname{cup}(\tilde G_{11,3})$, $\operatorname{cup}(\tilde G_{12,3})$ is
at least $5$.
\item"(d)" If $k\geq 4$, then
$$\operatorname{cup}(\tilde G_{n,k})\geq \cases
\dsize\frac{n-k+6}{2} &\text{ if } n-k+3\geq 7 \text{ is odd, } n-k+3\notin \{9,11\}, \\
\dsize\frac{n-k+5}{2} &\text{ if } n-k+3\geq 8 \text{ is even, } n-k+3\notin \{10,12\}, \\
5 &\text{ if } n-k+3\in \{9, 10, 11, 12\} .\\
\endcases $$

\item"(e)" If $k\geq 4$, then
$$\operatorname{cup}_{\Bbb Q}(\tilde G_{n,k})\geq \cases
1+\left[\dsize\frac{k}{2}\right]\left[\dsize\frac{n-k}{2}\right]
&\text{ if } 4\left[\dsize\frac{k}{2}\right]\left[\dsize\frac{n-k}{2}\right]<k(n-k), \\
\left[\dsize\frac{k}{2}\right]\left[\dsize\frac{n-k}{2}\right] &\text{ if
}4\left[\dsize\frac{k}{2}\right]\left[\dsize\frac{n-k}{2}\right]
=k(n-k).
\endcases $$
\endroster
\endproclaim
The author acknowledges that Proposition B(e) was suggested by Parames\-wa\-ran Sankaran.
Proposition B obviously implies the following (the upper estimates are the Grossman upper
bounds).

\proclaim{Corollary C} For the oriented Grassmann manifolds $\tilde G_{n,k}$, $n\geq 2k\geq
6$, we have: \roster
\item"(a)" $4\leq \operatorname{cat}(\tilde G_{6,3})\leq 5$;
\item"(b)" $\dsize\frac{n+5}{2}\leq \operatorname{cat}(\tilde G_{n,3})
\leq \dsize\frac{3n-7}{2}$ if $n$ is odd, $n\notin\{9, 11\}$, and
$\dsize\frac{n+4}{2}\leq \operatorname{cat}(\tilde G_{n,3}) \leq
\dsize\frac{3n-7}{2}$ if $n$ is even, $n\notin\{6, 10, 12\}$. In
addition to this, we have $6\leq \operatorname{cat}(\tilde
G_{9,3})\leq 10,\,\, 6\leq \operatorname{cat}(\tilde G_{10,3})\leq
11$, $6 \leq \operatorname{cat}(\tilde G_{11,3})\leq 13, \,\, 6
\leq \operatorname{cat}(\tilde G_{12,3})\leq 14$.
\item"(c)" If $k\geq 4$, then
$$\operatorname{cat}(\tilde G_{n,k})\geq \cases 6 &\text{ if } (n,k)=(8,4),\\
2+\left[\dsize\frac{k}{2}\right]\left[\dsize\frac{n-k}{2}\right] &\text{ if }
4\left[\dsize\frac{k}{2}\right]\left[\dsize\frac{n-k}{2}\right]<k(n-k), \\
1+\left[\dsize\frac{k}{2}\right]\left[\dsize\frac{n-k}{2}\right] &\text{ if } (n,k)\neq (8,4)
\text{ and } 4\left[\dsize\frac{k}{2}\right]\left[\dsize\frac{n-k}{2}\right] =k(n-k).\endcases $$

\endroster
\endproclaim

In Section 4, we shall prove the following upper bounds for $\operatorname{cup}(\tilde
G_{n,k})$ ($n\neq 6$)  and $\operatorname{cup}_{\Bbb Q}(\tilde G_{n,k})$ ($k\geq 4$); the
proofs can be seen as a realization of the third step of our Research Program. The manifold
$\tilde G_{6,3}$ will be treated in a special way; note that for this manifold the upper
bound coincides with the lower bound, so that we obtain $\operatorname{cup}(\tilde
G_{6,3})=3$. Also note that the upper bound for $\operatorname{cup}_{\Bbb Q}(\tilde G_{n,k})$
($k\geq 4$) given below coincides with the lower bound given in Proposition B(e) in
infinitely many cases.

\proclaim{Proposition D} For the oriented Grassmann manifolds $\tilde
G_{n,k}$, $n\geq 2k\geq 6$, we have

\roster
\item"(a)" $\operatorname{cup}(\tilde G_{6,3})\leq 3$. As a
consequence of this and Proposition {\text{\rm B(a)}, we have $\operatorname{cup}(\tilde
G_{6,3})=3$.

\item"(b)" We have
$$\align
&\operatorname{cup}(\tilde G_{n,3})\leq
\cases
\left[\dsize\frac{2^{s+2}-7}{3}\right] &\text{ if } n=2^s +1, s\geq 3,\\
\left[\dsize\frac{2^{s+2}-3}{3}\right] &\text{ if } n=2^s +2, s\geq 3,\\
\left[\dsize\frac{2^{s+2}+5\cdot 2^p-8}{3}\right] &\text{ if }
n=2^s+2^p +1, s>p\geq 1,\\
\left[\dsize\frac{2^{s+2}+5\cdot 2^p +3t-7}{3}\right] &\text{ if }
n=2^s+2^p+t+1, s>p\geq 1,\\&\,\,1\leq t\leq 2^p -1;
\endcases\endalign $$
$$
\align
&\operatorname{cup}(\tilde G_{n,4})\leq \cases
\left[\dsize\frac{5\cdot 2^{s}-13}{3}\right] &\text{ if } n=2^s +1, s\geq 3,\\
{2^{s+1}-4} &\text{ if } n=2^s +2, s\geq 3,\\
{2^{s+1}-3} &\text{ if } n=2^s +3, s\geq 3,\\
\left[\dsize\frac{2^{s+1}+4n-17}{3}\right] &\text{ if } 2^s+4\leq n\leq 2^{s+1};\\
\endcases\cr
&\operatorname{cup}(\tilde G_{n,k})\leq
\cases
\left[\dsize\frac{(k+1)\cdot 2^{s}+k-k^2-1}{3}\right] &\text{ if } n=2^s +1, s\geq 3, k\geq 5\\
\left[\dsize\frac{2^{s+1}+kn-k^2-1}{3}\right] &\text{ if } 2^s+2\leq n\leq 2^{s+1}, k\geq 5;\\
\endcases
\endalign
$$
\item"(c)" $\operatorname{cup}_{\Bbb Q}(\tilde G_{n,k})\leq \dsize \frac{k(n-k)}{4}$ if
$k\geq 4$. As a consequence of this and Proposition {\text{\rm B(e)}, we have
$\operatorname{cup}_{\Bbb Q}(\tilde G_{n,k})= \dsize\frac{k(n-k)}{4}$ for $n$ even and $k$
$(\geq 4)$ even, as well as for $n=4t+9$ $(t\geq 1)$ and $k=4$.
\endroster
\endproclaim
Note that Proposition B(d) gives $\operatorname{cup}(\tilde G_{8,4})\geq 5$, while
Proposition B(e) yields $\operatorname{cup}_{\Bbb Q}(\tilde G_{8,4})\geq 4$ and, from
Proposition D(c), we can see that $\operatorname{cup}_{\Bbb Q}(\tilde G_{8,4})=4$. This
indicates that, for a given $\tilde G_{n,k}$, the value of $\operatorname{cup}(\tilde
G_{n,k})$ may be higher than $\operatorname{cup}_{\Bbb Q}(\tilde G_{n,k})$. But the lower
bounds presently known for the $\Bbb Z_2$-cup-length of $\tilde G_{n,k}$ with $k\geq 4$ and
$(n,k)\neq (8,4)$, given in Proposition B(d), do not exceed the corresponding lower bounds
for the rational cup-length, given in Proposition B(e).

We hope that our Research Program based on Theorem A can lead to further interesting results
on the cup-length and L-S category not only for the oriented Grassmann manifolds but also for
other manifolds.

\head 2. ON THE CUP-LENGTH OF $R$-POINCAR\'E SPACES / PROOF OF THE MAIN RESULT
\endhead
In the spirit of W. Browder's \cite{4}, by an $R$-Poincar\'e space of formal dimension $n$ we
understand a path connected topological space $X$ for which there is an element $\mu\in
H_n(X)\cong R$ such that the $\cap$-product homomorphism $\cap\mu: H^k(X; R)\rightarrow
H_{n-k}(X; R)$, $x\mapsto x\cap \mu$, is an isomorphism for each $k$. By saying that an
$R$-Poincar\'e space $X$ is $R$-homologically $t$-connected ($t\geq 0$) we understand that
its reduced cohomology groups $\tilde H^i(X; R)$ vanish for all $i\leq t$. For example, any
$t$-connected (in the standard sense) closed $n$-dimensional manifold orientable over $R$ is
an $R$-homologically $t$-connected $R$-Poincar\'e space of formal dimension $n$.

\subhead 2.1 Proof of Theorem A\endsubhead If $X$ is an $R$-Poincar\'e space of formal
dimension $n$, then the cup product pairing $H^k(X; R)\times H^{n-k}(X; R) \rightarrow R$ is
nonsingular (see e.g. \cite{4} or \cite{9}); as a consequence, for any nonzero $x\in H^k(X;
R)$ there exists some $y\in H^{n-k}(X; R)$ such that $x\cup y$ is nonzero in $H^n(X; R)$. In
particular, this immediately implies (a3).

 For the rest of the proof, first note that the hypothesis of Theorem A implies that the
space $X$ is $R$-homologically $(r-1)$-connected and it has at least three nontrivial
unreduced $R$-cohomology groups. If there are just three, then Proposition A(a) is verified
in an obvious way. Indeed, in such a case we have $q=r=\dsize\frac{n}{2}$, and one readily
sees that $\operatorname{cup}_R(X)=2$.

So suppose now that $X$ has at least four
nontrivial unreduced $R$-cohomology groups, so that we have $q>r$.
Of course (see the beginning of this proof), any cup product
of the maximum length, that is, of the length
$\operatorname{cup}_R(X)$, must belong to $H^n(X; R)\cong R$.
So the cup-length of $X$ is realized by a cup product
$${x(1)_r^{p_1}}\cdots{x(s)_r^{p_s}}
{y_q^v}z_{q+i_1}^{j_1}\cdots z_{q+{i_m}}^{j_m}\in
H^n(X; R)\setminus \{0\}, \tag $\bullet$ $$ where
$x(1)_r, \dots, x(s)_r\in H^r(X; R)$, $y_q\in H^q(X; R)$, $z_{q+{i_l}}\in
H^{q+{i_l}}(X; R)$ are nonzero cohomology classes and $p_1,\dots, p_s$,
$v$, $j_1,\dots,j_m$, $i_1, \dots, i_m$ are non-negative integers.
Denote $p=p_1+\dots+p_s$. Then, of course,
$\operatorname{cup}_R(X)=p+v+{j_1}+\dots+{j_m}$.

From this it is clear that
$$pr+vq+{j_1}(q+i_1)+\dots+{j_m}(q+i_m)=n,$$
therefore
$$n\geq r(p+v+{j_1}+\dots+{j_m}).$$ In other words, we obtain
$$n\geq r\cdot \operatorname{cup}_R(X),$$
which proves (a1). If there is a cohomology class $x\in H^r(X; R)\setminus \{0\}$ such that
$r\cdot\operatorname{ht}_R(x)=n$, then obviously $\operatorname{cup}_R(X)\geq
\dsize\frac{n}{r}$; this together with (a1) proves (a2). The rest of Theorem A(a) is clear in
view of what we have said in the beginning of this proof.

We pass to the proof of Theorem A(b). Now we suppose that $r<q$, and we have a basis $a(1)_r,\dots,
a(t)_r$ for $H^r(X; R)$ such that
$$a(1)_r^{k_1+1}=0,\dots, a(t)_r^{k_t+1}=0$$
for some positive integers $k_1,\dots, k_t$ such that
$r(k_1+\dots+k_t)<n$. As above, take an element realizing the
cup-length of $X$, hence some
$$c:={x(1)_r^{{\pi}_1}}\cdots{x(s)_r^{{\pi}_s}}
{y_q^v}z_{q+i_1}^{j_1}\cdots z_{q+{i_m}}^{j_m}\in
H^n(X; R)\setminus \{0\}, \tag $\bullet$ $$ where
$x(1)_r, \dots, x(s)_r\in
H^r(X; R)\setminus \{0\}$, $y_q\in H^q(X; R)\setminus \{0\}$,
$z_{q+{i_l}}\in H^{q+{i_l}}(X; R)\setminus
\{0\}$, and ${\pi}_1,\dots, {\pi}_s$, $v$, $j_1,\dots,j_m$, $i_1, \dots,
i_m$ are non-negative integers. We denote $\pi=\pi_1+\dots+\pi_s$. But now
$x(i)_r=\sum_{j=1}^t{\alpha_{i,j}} a(j)_r$ for some uniquely determined
coefficients $\alpha_{i,1},\dots,\alpha_{i,t}\in R$. So
$c$ is a linear combination of cup products of the form
$$a(1)_r^{p_1}\cdots a(t)_r^{p_t}{y_q^v}z_{q+i_1}^{j_1}\cdots
z_{q+{i_m}}^{j_m}, \tag $\bullet\bullet$ $$ where $p_1,\dots,p_t$ are non-negative integers
such that $p_1+\dots+p_t=\pi$.
Since $c\neq 0$, at least one of the products $(\bullet\bullet)$
must be nonzero; of course, in such a nonzero cup product, the exponent of $a(i)_r$ must be
less than or equal to $k_i$ for all $i=1,\dots,t$. We conclude that the cup-length is
realized by a nonzero element
$$a(1)_r^{p_1}\cdots a(t)_r^{p_t}{y_q^v}z_{q+i_1}^{j_1}\cdots
z_{q+{i_m}}^{j_m}\in H^n(X; R)\setminus \{0\}, \tag $\bullet\bullet\bullet$ $$ where
$p_1,\dots,p_t, v, {j_1},\dots,{j_m}, {i_1},\dots,{i_m}$ are non-negative integers, $p_1\leq
k_1,\dots, p_t\leq k_t$, $\operatorname{cup}_R(X)=p_1+\dots+p_t+v+{j_1}+\dots+{j_m}$, and
$y_q\in H^q(X; R)\setminus \{0\}$, $z_{q+{i_l}}\in H^{q+{i_l}}(X; R)\setminus \{0\}$. Using
the fact that $q-r>0$, we obtain
$$\eqalign{n&=r(p_1+\dots+p_t)+q(v+j_1+\dots+j_m)+{i_1}{j_1}+\dots+{i_m}{j_m}\cr
&=q(p_1+\dots+p_t+v+j_1+\dots+j_m)+(r-q)(p_1+\dots+p_t)+{i_1}{j_1}+\dots+{i_m}{j_m}\cr &\geq
q(p_1+\dots+p_t+v+j_1+\dots+j_m)+(r-q)(p_1+\dots+p_t)\cr
&=q\operatorname{cup}_R(X)-(q-r)(p_1+\dots+p_t)\cr &\geq
q\operatorname{cup}_R(X)-(q-r)(k_1+\dots+k_t).}$$ The proof of Theorem A is finished.

\head 3. APPLICATIONS: LOWER BOUNDS FOR $\tilde G_{n,k}$ / PROOF OF PROPOSITION B \endhead

\subhead 3.1 Selected facts about the oriented Grassmann
manifolds\endsubhead
 The oriented Grassmann manifold $\tilde G_{n,k}$ consists of {\it oriented} $k$-dimensional
 vector subspaces in $\Bbb R^n$, similarly the Grassmann manifold $G_{n,k}$ consists of
{\it all} $k$-dimensional vector subspaces in $\Bbb R^n$. For obvious reasons, we shall
suppose that $2k\leq n$ and, for the reason explained in the Introduction, we restrict
ourselves to $k\geq 3$.

Let $\gamma_{n,k}$ (briefly $\gamma$) denote the canonical $k$-plane bundle over $G_{n,k}$,
$\tilde \gamma_{n,k}$ (briefly $\tilde \gamma$) be the canonical $k$-plane bundle over
$\tilde G_{n,k}$, and $\xi$ be the canonical line bundle associated with the double covering
$p: \tilde G_{n,k} \rightarrow G_{n,k}$. We denote $w_i=w_i(\gamma)$ in $H^i(G_{n,k}; \Bbb
Z_2)$ and $\tilde w_i=w_i(\tilde\gamma)$ in $H^i(\tilde G_{n,k}; \Bbb Z_2)$ the corresponding
Stiefel-Whitney classes. Of course, $\tilde w_1=0$, and $w_1$ is easily seen to coincide with
the first Stiefel-Whitney class of $\xi$.

By Borel \cite{2}, the $\Bbb Z_2$-cohomology ring $H^\ast(G_{n,k}; \Bbb Z_2)$ can be
identified with a quotient ring,
$$H^\ast(G_{n,k}; \Bbb Z_2)\cong \Bbb Z_2[w_1,\dots, w_k]/I_{n,k},$$
where the ideal $I_{n,k}$ is generated by the homogeneous
components of
$$\frac{1}{1+w_1+\dots+w_k}$$
in dimensions $n-k+1,\dots,n$.

Calculations in the ring $H^\ast(\tilde G_{n,k}; \Bbb Z_2)$ ($k\geq 3$) may be very
difficult. Fortunately, they can be made a little easier thanks to the Gysin exact sequence
associated with the double covering $p: \tilde G_{n,k} \rightarrow G_{n,k}$. Indeed, for
showing that $\tilde w_2^{i_2}\cdots \tilde w_k^{i_k}\neq 0$, it is enough to verify that
$w_2^{i_2}\cdots w_k^{i_k}\in H^{2i_2+\dots+ki_k}(G_{n,k}; \Bbb Z_2)$ is not a multiple of
the class $w_1$.

In view of Theorem A (and in view of the step (3) of our Research Program), it is good to
keep in mind an explicit description of the first nontrivial reduced $\Bbb Z_2$-cohomology
group for those oriented Grassmann manifolds we are interested in.

\proclaim{Lemma E} For the oriented Grassmann manifolds $\tilde G_{n,k}$ ($n\geq 2k\geq 6$),
we have $H^2(\tilde G_{n,k}; \Bbb Z_2) =\{0, \tilde w_2\}\cong \Bbb Z_2$.
\endproclaim
\demo{Proof} Let $f: V_{n,k}\rightarrow \tilde G_{n,k}$ (where $V_{n,k}$ is the Stiefel manifold of
orthonormal $k$-frames in $\Bbb R^n$) be the natural fiber bundle with fiber $SO(k)$. The Lemma is
readily implied by the corresponding Serre exact sequence.
\enddemo

\subhead 3.2 Preparations for the proof of Propositions B(a)-(c)
\endsubhead To any integer $n\geq 6$, we find $s$ as the unique integer such
that $2^s<n\leq 2^{s+1}$. Since $\operatorname{dim}(G_{n,3})<4n\leq 2^{s+3}$, we have
$$(1+w_2+w_3)^{2^{s+3}}=1. \tag *$$
To show that
$$w_2^{i_2}w_3^{i_3}\in H^{2i_2+3i_3}(G_{n,3}; \Bbb Z_2)$$ cannot be expressed
as a multiple of the class $w_1$, it is enough to show that $w_2^{i_2}w_3^{i_3}$ is not in
the ideal $J_{n,3}$ of $\Bbb Z_2[w_2,w_3]$ generated by the homogeneous components of
$$\frac{1}{1+w_2+w_3}=(1+w_2+w_3)^{2^{s+3}-1}$$
(we have used (*)) in dimensions $n-2$, $n-1$, $n$.
 Therefore the ideal $J_{n,3}$ is generated by homogeneous polynomials $g_{n-2}$, $g_{n-1}$,
$g_{n}$, obtained as the homogeneous components in dimensions $n-2$, $n-1$ and $n$,
respectively, of
$$
{\sum\limits_{i=0}^{2^{s+3}-1}}{\sum\limits_{j=0}^{i}}\,\,{{i}\choose{j}} w_2^{i-j}w_3^j. \tag **
$$
Note that the binomial coefficient ${2^{s+3}-1}\choose{i}$ is odd for every
$i=0,1,\dots,2^{s+3}-1$. For $\kappa=n-2,\, n-1,\, n$ one calculates from (**) that
$$g_{\kappa}=\sum\limits_{\frac{\kappa}{3}\leq i\leq \frac{\kappa}{2}}{
{i}\choose{3i-\kappa}} w_2^{3i-\kappa}w_3^{\kappa-2i}. \tag *** $$

\subhead 3.3 Proof of Proposition B(a)\endsubhead From (***) we obtain that the ideal
$J_{6,3}$ in $\Bbb Z_2[w_2,w_3]$ is generated by $g_4=w_2^2$  and $g_6=w_3^2+w_2^3$. So $w_2
w_3$ is not in $J_{6,3}$. Consequently, $\tilde w_2\tilde w_3\neq 0$, and (cf. (a3) in
Theorem A) $\operatorname{cup}(\tilde G_{6,3})\geq 3$, as we have claimed.

\subhead 3.4 Proof of Proposition B(b)\endsubhead We say that a homogeneous polynomial
$h_a\in \Bbb Z_2[w_2,w_3]$, in dimension $a$, is $w_2$-monomial (\lq\lq monomial\rq\rq\,\, is
an adjective here) if $h_a=w_2^{\frac{a}{2}}$. As a realization of the step (2) of our
Research Program, we should find $a$ as high as we can (this value of $a$ will be called an
{\it available target dimension}) such that no element in the $a$-dimensional homogeneous
component of $J_{n,3}$ is $w_2$-monomial. This then gives that $\tilde w_2^{\frac{a}{2}}\neq
0$, and Proposition B(b) is readily implied by (a2) or (a3) of Theorem A.

{\it Odd values of} $n$. Let us first suppose that $n$ is odd, $n\geq 7$, $n\notin \{9,11\}$.
We shall show that $n+1$ is an available target dimension; so our aim now is to verify that
none of the elements $w_3g_{n-2}$, $w_2g_{n-1}$, $w_3g_{n-2}+w_2g_{n-1}$ is $w_2$-monomial.

If $n=6t +1$, then using (***) we calculate that
$$w_2g_{n-1}=w_2w_3^{2t}+\dots+w_2^{3t+1},$$
$$w_3g_{n-2}=0\cdot w_2w_3^{2t}+\dots+
(3t-1)w_2^{3t-2}w_3^2,$$  and
$$w_2g_{n-1}+w_3g_{n-2}= w_2w_3^{2t}+\dots+w_2^{3t+1}.$$
Hence none of $w_3g_{n-2}$, $w_2g_{n-1}$, $w_3g_{n-2}+w_2g_{n-1}$
is $w_2$-monomial, indeed.

If $n=6t +3$, we shall suppose that $t>1$; the case $n=9$ has a separate treatment. Write
$n=3\cdot 2^{k+1}(2l+1)+3$. Using (***), one readily verifies that none of the polynomials
$w_2g_{n-1}$, $w_3g_{n-2}$ is $w_2$-monomial. Further, one calculates that
$$w_2g_{n-1}+w_3g_{n-2}=
\sum_{i=2t+1}^{3t} {{i+1}\choose {3i-6t-1}}
w_2^{3i-6t-1}w_3^{6t+2-2i} + 1\cdot w_2^{3t+2}.
$$
Therefore, \roster
\item"(a)" if $k\geq 1$ and $l$ is arbitrary, then
$$w_2g_{n-1}+w_3g_{n-2}=\dots+w_2^2w_3^{2^{k+1}(2l+1)}+\dots+
w_2^{3\cdot 2^k(2l+1)+2},$$

\item"(b)" if $k=0$ and $l$ is at least $2$ and even, then
$$w_2g_{n-1}+w_3g_{n-2}=\dots+w_2^5w_3^{4l}+\dots+w_2^{6l+5},$$

\item"(c)" if $k=0$ and $l$ is at least $1$ and odd, then
$$w_2g_{n-1}+w_3g_{n-2}=\dots+w_2^{6l+2}w_3^{2}+w_2^{6l+5}.$$
\endroster
We conclude that for $n=6t+3$ with $t>1$, the polynomial $w_2g_{n-1}+w_3g_{n-2}$ is not
$w_2$-monomial.

For $n$ odd, we are left with the case $n=6t+5$; we suppose
that $t>1$ ($n=11$ is treated separately).
We write now $n=3\cdot 2^{k+1}\cdot (2l+1)-1$.

If $l\geq 1$, then we obtain from (***) that
$$w_2g_{n-1}=w_2^{3\cdot 2^k}w_3^{3\cdot 2^{k+2}\cdot l}+
\dots + w_2^{\frac{n+1}{2}},$$
$$w_3g_{n-2}=w_3^{\frac{n+1}{3}}+\dots+
(3\cdot 2^k(2l+1)-2)w_2^{3\cdot 2^k(2l+1)-3}w_3^2.$$ Hence for $l\geq 1$ none of the polynomials
$w_2g_{n-1}$, $w_3g_{n-2}$, $w_2g_{n-1}+w_3g_{n-2}=w_3^\frac{n+1}{3}+\dots+w_2^{\frac{n+1}{2}}$ is
$w_2$-monomial.

We are left with $n=2^{k+1}\cdot 3\cdot (2l+1)-1$ for $l=0$; now
the argument used for $l\geq 1$ does not work. Note
that we take $k>1$, because we have $n >11$.

Now one calculates from (***) that
$$g_{n-2}=w_3^{\frac{n-2}{3}}+\dots+0\cdot w_2^{\frac{n-5}{2}}w_3.$$
Further, if $k$ is odd, then
$$w_2g_{n-1}=\dots+w_2^{\frac{n+7}{3}}w_3^{\frac{n-11}{9}}+
\dots +w_2^{\frac{n+1}{2}},$$
and if $k$ is even, then
$$w_2g_{n-1}=\dots+w_2^{\frac{n+4}{3}}w_3^{\frac{n-5}{9}}+
\dots +w_2^{\frac{n+1}{2}}.$$ We conclude that none of the polynomials $w_2g_{n-1}$,
$w_3g_{n-2}$, $w_2g_{n-1}+w_3g_{n-2}= w_3^{\frac{n+1}{3}}+\dots+ w_2^{\frac{n+1}{2}}$ is
$w_2$-monomial.

{\it Even values of} $n$. Let us suppose that $n$ is even, $n\geq 8$, $n\notin \{10, 12\}$.
By a similar analysis as above, we could show that $n$ is an available target dimension. But
the hard work can be avoided. Indeed, we have the standard \lq\lq inclusion''\,\, $i:\tilde
G_{n-1,3}\rightarrow \tilde G_{n,3}$ such that
$i^\ast(\tilde\gamma_{n,3})=\tilde\gamma_{n-1,3}$. Since $n-1$ is odd, we know, by what we
have computed above, that $n$ is an available target dimension for $\tilde G_{n-1,3}$. This
implies that $\tilde w_2^{\frac{n}{2}}\in H^n(\tilde G_{n-1,3}; \Bbb Z_2)$ does not vanish,
and therefore $\tilde w_2^{\frac{n}{2}}\in H^n(\tilde G_{n,3}; \Bbb Z_2)$ cannot be zero.
Since $n$ is now less than $\operatorname{dim}(\tilde G_{n,3})$, (a3) of Theorem A applies.
Proposition B(b) is proved.

\subhead 3.5 Proof of Proposition B(c)\endsubhead
 To finish the proof of Proposition B, we are left with the special cases,
$n=9,\, 10,\, 11,\, 12$. For $n=9$, we obtain that $g_7=w_2^2w_3$, $g_8=w_2w_3^2+w_2^4$,
$g_9=w_3^3$. So $\tilde w_2^4\neq 0$, while $\tilde w_2^5=0$ (indeed: we have
$w_2^5=w_3g_7+w_2g_8$). In the three remaining cases, it is then readily seen that $\tilde
w_2^4$ is not zero (apply the \lq\lq inclusions\rq\rq\,\, $i:\tilde G_{9,3}\rightarrow \tilde
G_{9+s,3}$ such that $i^\ast(\tilde\gamma_{9+s,3})=\tilde\gamma_{9,3}$, $s=1,\,2,\,3$).
Proposition B(c) is proved.

\subhead 3.6 Proof of Proposition B(d)\endsubhead Let $j: \tilde G_{a,b}\rightarrow \tilde
G_{a+1,b+1}$ be the \lq\lq inclusion\rq\rq\,\, such that $j^\ast(\tilde \gamma_{a+1,b+1}) =
\tilde\gamma_{a,b}\oplus\varepsilon^1$, where $\varepsilon^1$ is the trivial line bundle.
Then we have $j^\ast(w_2(\tilde \gamma_{a+1,b+1}))=w_2(\tilde \gamma_{a,b})$. By an obvious
iteration, we obtain the \lq\lq inclusion\rq\rq\,\, $i:\tilde G_{n-k+3,3}\rightarrow \tilde
G_{n,k}$ such that $i^\ast(w_2(\tilde \gamma_{n,k}))=w_2(\tilde\gamma_{n-k+3,3})$. In the
proofs of Propositions B(b), B(c), we have shown that certain powers of $w_2(\tilde
\gamma_{n-k+3,3})$ do not vanish; of course, then the same powers of $w_2(\tilde
\gamma_{n,k})$ do not vanish, and the rest is clear: one applies (a3) from Theorem A.
Proposition B(d) is proved.

\subhead 3.7 Proof of Proposition B(e)\endsubhead Now $k\geq 4$ and, using the known
description of the cohomology algebra $H^\ast(\tilde G_{n,k}; \Bbb Q)$ (cf. for instance
\cite{11; Proposition 5}), one readily verifies that we have $r=4$ in Theorem A. Let
$p_i(\tilde \gamma_{n,k})\in H^{4i}(\tilde G_{n,k}; \Bbb Q)$ be the $i$th rational Pontrjagin
class. The height of $p_1(\tilde \gamma_{n,k})$ is
$\left[\frac{k}{2}\right]\left[\frac{n-k}{2}\right]$ (see e.g. the proof of \cite{11; Theorem
1}). The lower bounds stated in Proposition B(e) are then implied by (a2) or (a3) of Theorem
A.

The proof of Proposition B is finished.

\head 4. APPLICATIONS: UPPER BOUNDS FOR $\tilde G_{n,k}$ / PROOF OF PROPOSITION D
\endhead \subhead 4.1 Proof of Proposition D(a)\endsubhead From the proof of Proposition B(a),
we know that $\tilde w_2\tilde w_3\neq 0$. Hence there exists a nonzero cohomology class
$a\in H^4(\tilde G_{6,3}; \Bbb Z_2)$ such that $\tilde w_2\cdot \tilde w_3\cdot a\neq 0$. It
is clear that $\operatorname{cup}(\tilde G_{6,3})\geq 3$, and that $\operatorname{cup}(\tilde
G_{6,3})$ could be more than $3$ only if $a$ could be decomposed as a product of two
cohomology classes (in view of Lemma E in 3.1, the only decomposition, not excluded a priori,
would be $a=\tilde w_2^2$) or if $\tilde w_2^4$ would be nonzero. But the element $a$ is
indecomposable, because, in addition to the fact that $H^1(\tilde G_{6,3}; \Bbb Z_2)=0$, we
have (Lemma E) $H^2(\tilde G_{6,3}; \Bbb Z_2)=\{0, \tilde w_2\}$, and $\tilde w_2^2=0$ (cf.
3.3). So we conclude that $\operatorname{cup}(\tilde G_{6,3})\leq 3$. Proposition D(a) is
proved.

\subhead 4.2 Proof of Proposition D(b)\endsubhead We apply Theorem A(b), with $r=2$, $R=\Bbb
Z_2$. In view of Lemma E, we have $\tilde w_2$ as the only choice of basis in $H^2(\tilde
G_{n,k}; \Bbb Z_2)$ ($n\geq 2k\geq 6$). Dutta and Khare \cite{6} calculated the height of
$w_2$ in $H^\ast(G_{n,k}; \Bbb Z_2)$ as follows (we quote just the results we need, hence for
$n\geq 2k\geq 6$).

\proclaim{Lemma F} {\rm(}S. Dutta, S. Khare \cite{6}{\rm)} For the Grassmann manifolds
$G_{n,k}$ ($n\geq 2k\geq 6$) one has
$$\align
&\operatorname{ht}(w_2(\gamma_{n,3}))=
\cases
{2^s}-1&\text{ if } n={2^s}+1,\\
2^s &\text{ if } n=2^s +2,\\
{2^s}+{2^{p+1}}-2&\text{ if }
n=2^s+2^p +1, s>p\geq 1,\\
{2^s}+{2^{p+1}}-1&\text{ if }
n=2^s+2^p+t+1, s>p\geq 1,\\&\,\,1\leq t\leq 2^p -1;
\endcases\cr
&\operatorname{ht}(w_2(\gamma_{n,4}))=
\cases
{2^s}-1&\text{ if } n=2^s +1,\\
{2^{s+1}-4} &\text{ if } n=2^s +2,\\
{2^{s+1}-4} &\text{ if } n=2^s +3,\\
{2^{s+1}}-1&\text{ if } 2^s+4\leq n\leq 2^{s+1};\\
\endcases\cr
&\operatorname{ht}(w_2(\gamma_{n,k}))=
\cases
{2^s}-1&\text{ if } n=2^s +1, k\geq 5\\
{2^{s+1}}-1&\text{ if } 2^s+2\leq n\leq 2^{s+1}, k\geq 5.\\
\endcases
\endalign
$$
\endproclaim

It is clear that $\operatorname{ht}(w_2(\gamma_{n,k}))$, briefly $\operatorname{ht}(w_2)$,
cannot exceed half of the dimension of the corresponding Grassmann manifold. But, as we see
from Lemma F, it is mostly smaller. In an obvious way, the above quoted results on
$\operatorname{ht}(w_2)$ enable us to find, for each pair $(n,k)$ under consideration, some
$c$, mostly smaller than half of the dimension, such that $\tilde w_2^{c+1}=0$. For instance,
for $G_{2^s+1,3}$ we know by Lemma F that $\operatorname{ht}(w_2)=2^s -1$. Therefore
$w_2^{2^s}=0$, and of course for $\tilde G_{2^s+1,3}$ we have $\tilde w_2^{2^s}=0$. Realizing
the step (3) of our Research Program, using Theorem A(b), we obtain the upper bound stated in
Proposition D(b) for this case. The remaining cases are similar: when
$\operatorname{ht}(w_2)$ is less than half of the dimension, then we obtain the upper bound
stated in Proposition D(b) by applying Theorem A(b), and in the cases where
$\operatorname{ht}(w_2)$ is precisely half of the dimension, we have just half of the
dimension as an upper bound for $\operatorname{cup}(\tilde G_{n,k})$ from Theorem A(a).

\subhead 4.3 Proof of Proposition D(c)\endsubhead The upper bound given in Proposition D(c)
is obtained from (a1) of Theorem A (with $r=4$).

Proposition D is proved.

\head 5. REMARKS \endhead

\subhead 5.1 Remark\endsubhead In a special case, for the real flag manifolds
$F(n_1+\dots+n_s)=O(n_1+\dots+n_s)/O(n_1)\times\dots\times O(n_s)$, we derived an upper bound
of the same type as (b1) of Theorem A, with $r=1$ and $R=\Bbb Z_2$, in \cite{10; Proposition
3.2.2}. We based it there on specific properties of the $\Bbb Z_2$-cohomology of the flag
manifolds, but we did not recognize the full potential, now expressed in Theorem A(b). In
\cite{10}, we also illustrated the strength of that special case of (b1) by calculating the
exact value of the $\Bbb Z_2$-cup-length of $F(1,2,n_3)$ for all $n_3\geq 3$.

\subhead 5.2 Remark\endsubhead As pointed out by Akira Kono after having seen an early
version of the author's calculations for $\tilde G_{n,3}$, one can deal with $G_2/SO(4)$ in a
similar way (note that $G_2/SO(4)$ also is a simply-connected irreducible compact Riemannian
symmetric space). Using Borel and Hirzebruch's \cite{3; 17.3}, one calculates that
$\operatorname{cup}(G_2/SO(4))=4$. As a consequence, $\operatorname{cat}(G_2/SO(4))\geq 5$.
On the other hand, the Grossman upper bound yields $\operatorname{cat}(G_2/SO(4))\leq 5$, and
so $\operatorname{cat}(G_2/SO(4))=5$.

\subhead 5.3 Remark\endsubhead In view of Theorem A(b), it would be interesting to know the
exact values of $\operatorname{ht}(\tilde w_2)$.

\subhead 5.4 Remark\endsubhead For the oriented Grassmann manifolds $\tilde G_{n,k}$ we know
that $2\operatorname{ht}(\tilde w_2)<k(n-k)$ whenever $n$ is odd, independently of
Dutta-Khare's Lemma F. Indeed, if $n$ is odd, then (see e.g. \cite{1; Theorem 1.1})
$w_2(\tilde G_{n,k})=\tilde w_2$. Hence the value of ${\tilde w_2}^{\frac{k(n-k)}{2}}$ on the
fundamental class of the manifold $\tilde G_{n,k}$ is one of its Stiefel-Whitney numbers.
But, as is well known, all the Stiefel-Whitney numbers of $\tilde G_{n,k}$ vanish.

\remark{Acknowledgement} The author thanks A. Kono, J. L\"orinc, T. Macko, M. Mimura, P.
Sankaran, and R. Stong for their comments which were helpful at various stages of working on
this paper.

\endremark

\enddocument